\documentclass[journal]{IEEEtran}

\usepackage{graphicx}
\usepackage{amsmath,amssymb}
\usepackage[official]{eurosym}
\usepackage{multirow}
\usepackage{multicol}

\usepackage{array}
\usepackage[noadjust]{cite}
\usepackage[justification=centering, font=scriptsize]{caption}
\usepackage{booktabs}

\begin{document}

\title{A Prosumer-Centric Framework for Concurrent Generation and Transmission Planning — Part II}

\author{Ni~Wang,
        Remco~Verzijlbergh,
        Petra~Heijnen,
        and~Paulien~Herder
\thanks{N. Wang, R. Verzijlbergh and P. Heijnen is with the Faculty of Technology, Policy and Management, Delft University of Technology, Delft, the Netherlands (e-mails: \{n.wang, r.a.verzijlbergh, p.w.heijnen\}@tudelft.nl). }
\thanks{P. Herder is with the Faculty of Applied Sciences, Delft University of Technology, Delft, the Netherlands (e-mail: p.m.herder@tudelft.nl). }
\thanks{This research received funding from the Netherlands Organisation for Scientific Research (NWO) [project number: 647.002.007].}

}


\maketitle

\begin{abstract}
We propose a framework where generation and transmission capacities are planned concurrently in market environments with a focus on the prosumers. 
This paper is a continuation of Part I and presents numerical results from three archetypal case studies. 
Following the proposed framework, optimal planning decisions are shown in all the cases. Furthermore, in case study I, no-regret planning decisions considering the uncertainties in future electricity market designs are discussed. In case study II, we look at a situation where the social resistance of wind energy prevails and the prosumers choose to not invest in wind energy. This preference results in an increased system cost, and thus may harm other prosumers. The framework is used as a negotiation simulator to deal with this problem. Case study III presents numerical results for a mixed bilateral/pool market. The case studies utilize realistic data from the Dutch power system and the European power system to provide policy-relevant results that aid their decarbonization in various market environments.
\end{abstract}


%
\IEEEpeerreviewmaketitle

%
%
%
%

\section{Introduction}

This paper is a continuation of its Part I. There, we provide a prosumer-centric approach for concurrent generation and transmission planning, where 1) the costs and benefits for the prosumers along the planning horizon are explicitly modeled, 2) the investment preferences of the prosumers are included, 3) distributed algorithms are provided to solve the planning problem in a decentralized manner. The formulations of the investment problems integrating three different electricity markets, i.e. the P2P market, the pool market and the mixed bilateral/pool market, respectively, are given. 

We position this work as a general modeling framework to systematically assess the influences of various electricity market designs on the planning decisions. Therefore, the basic usage of our framework is to investigate the optimal investment decisions and the associated costs and benefits given the different market environments. By this, the framework is able to provide recommendations on issues that pertain to social, political and policy aspects.
To further illustrate its applicability, in this paper, three archetypal case studies will be presented.
These case studies focus on different aspects of the framework and will be done in various realistic contexts. Now they will be briefly introduced and more details will be given in the corresponding sections.

Case study I discusses the no-regret planning decisions taking into account the uncertainty in electricity market designs of the future. Case study II focuses on the effects of the investment preferences of the prosumers. There, the framework is used as a simulator for the negotiation that is incurred by the conflicting interests of the prosumers. Case study III presents numerical results for a mixed bilateral/pool market. The first two archetypal case studies are based on the national RES planning process in the Netherlands for 2030, while the last one takes a European perspective and models a highly renewable Europe in 2050.   

In Section \ref{sec:case_I} - \ref{sec:case_III}, the three case studies are described and the numerical results are presented, respectively. Section \ref{sec:conclusion} concludes the paper.

\section{Case Study I} \label{sec:case_I}

To combat climate change, the Dutch government aims to reduce the Netherlands’ greenhouse gas emissions by 49\% by 2030 compared to 1990 levels. In the electricity sector, a major focus is to increase the share of wind and solar energy, i.e. from 14.25 TWh (12\% of the total production) in 2018 to 84 TWh (72\% of the total production) in 2030. Accordingly, massive investments in wind and solar energy are needed. The investments are facilitated through a national program Regionale Energiestrategie (Regional Energy Strategies) \cite{UnievanWaterschappen2019}. In this program, the country is divided into 30 energy regions, where each region proposes its own investments in onshore wind and solar energy and meanwhile, the TSO proposes the transmission network expansion plan. 

This problem is dealt with using our framework. The generation and transmission planning problem is considered as optimization problems, where each agent minimizes their own cost concurrently. The cost for a region includes not only the costs that are considered in energy system planning models (i.e. CapEx, FOM and VOM), but also the expected trading costs in the electricity markets. 
In this case study, the influences of two market designs on the planning decisions will be investigated to account for the uncertainty in future market designs. In addition to a pool market, a P2P market as a promising next-generation electricity market will be discussed as the alternative future market design. 

We start by describing the model set-up in Section \ref{case_I_model}. Then, the optimal investment decisions and the relevant costs under different markets are presented in Section \ref{case_I_investment}. At last, in Section \ref{sec:case_I_uncertainty}, how the uncertainty is taken into account is described and the no-regret investment decisions are shown.

\subsection{Model set-up} \label{case_I_model}
The Netherlands is divided into 30 nodes which are connected by transmission lines. The considered time horizon is one year with hourly resolution. 
The input data includes, among others, the spatio-temporal variations of wind and solar capacity factors and their maximum potentials which is obtained from \cite{Wang2020b}. 
The considered technologies are onshore wind turbines, solar PV, OCGT, hydrogen storage and flow battery storage. Due to the differences in their efficiencies and cost parameters, hydrogen storage is best suited for long-term usage, whereas battery storage provides an option for short-term storage.
Furthermore, in this national planning process, the regions are not asked to propose capacity in offshore wind and thus offshore wind is not considered.  
Since quite a significant generation expansion is needed, for the ease of presenting results, the existing generation and storage capacities are not considered. The existing transmission network is considered and the capacity can be expanded. Direct current (DC) power flow calculations are performed based on the existing capacities. 
2030 estimations for the techno-economic parameters are used and are taken from \cite{Schlachtberger2017} except for the network cost parameter which is taken from \cite{Wang2020b}.
 
Regarding the product differentiation term in P2P markets, two points of view may be used \cite{Sorin2019}. The first one is the system view. The product differentiation can be regarded as an externality cost such as transaction cost (TC) or exogenous grid cost. The second view is the prosumer view. It represents the willingness to pay of the prosumers in relation to the trades. In this case study, we employ the first one, it will be used to represent the TC in the P2P markets. Later, in the second case study, the perspective of the latter will be investigated.

\subsection{Optimal investment decisions under different market designs} \label{case_I_investment}

We first look at the optimal investment decisions in terms of the total installed capacities for the system. Here, two different markets are considered, i.e. the pool market and the P2P market. We start with the pool market. Then, the P2P markets are analyzed with increasing TCs, starting from 10 \% of the average electricity price to 30 \%, 50 \%, 70 \% and 90 \% to investigate the influences of TC on the planning decisions.

\begin{figure}[htbp]
\centerline{\includegraphics[width=0.5\textwidth]{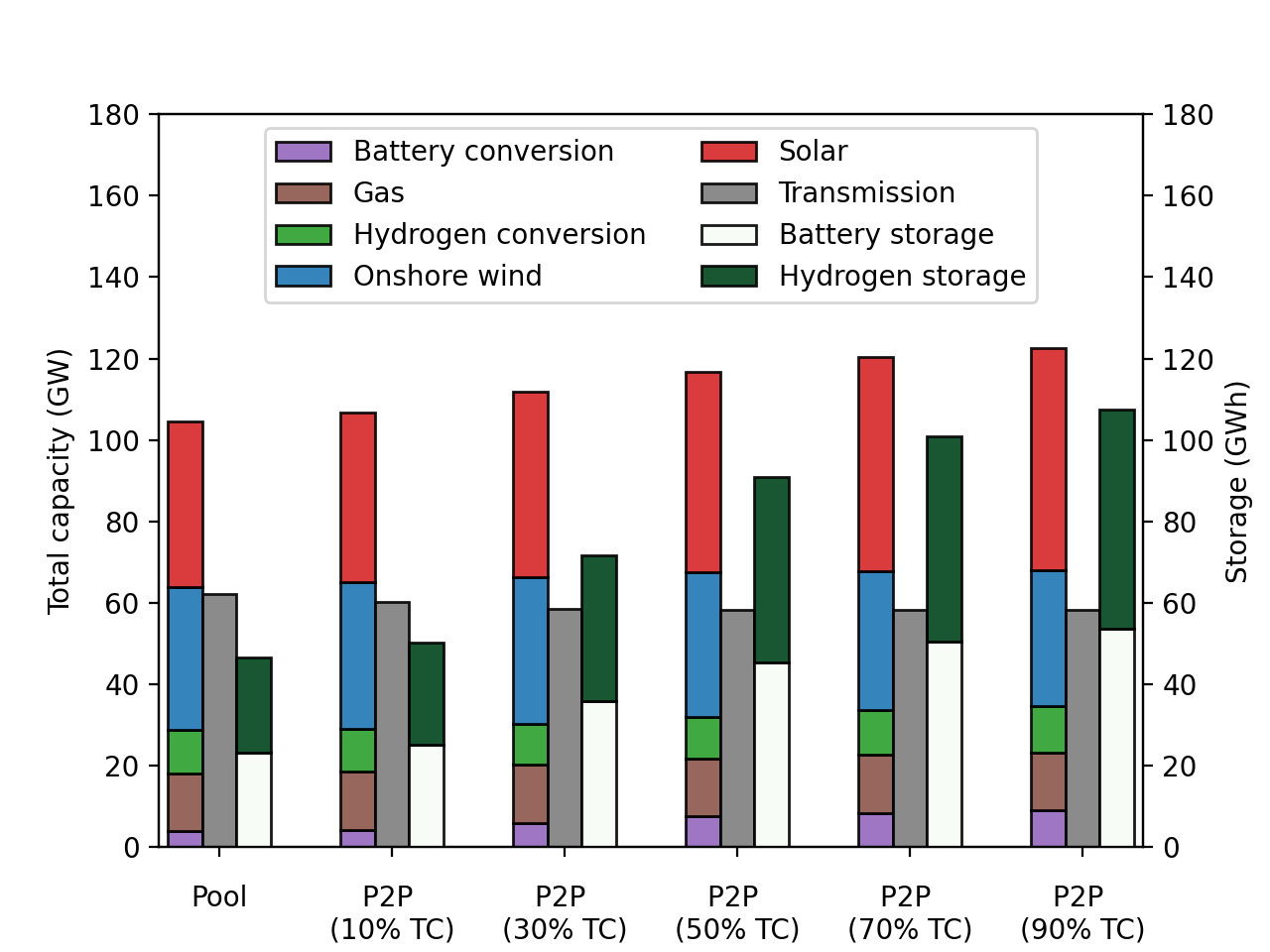}}
\caption{Optimal total capacities for the Netherlands under pool market and P2P market with different TC (\% of the average electricity price).}
\label{fig:case_I_total_capacity}
\end{figure}

\begin{figure}[htbp]
\centerline{\includegraphics[width=0.5\textwidth]{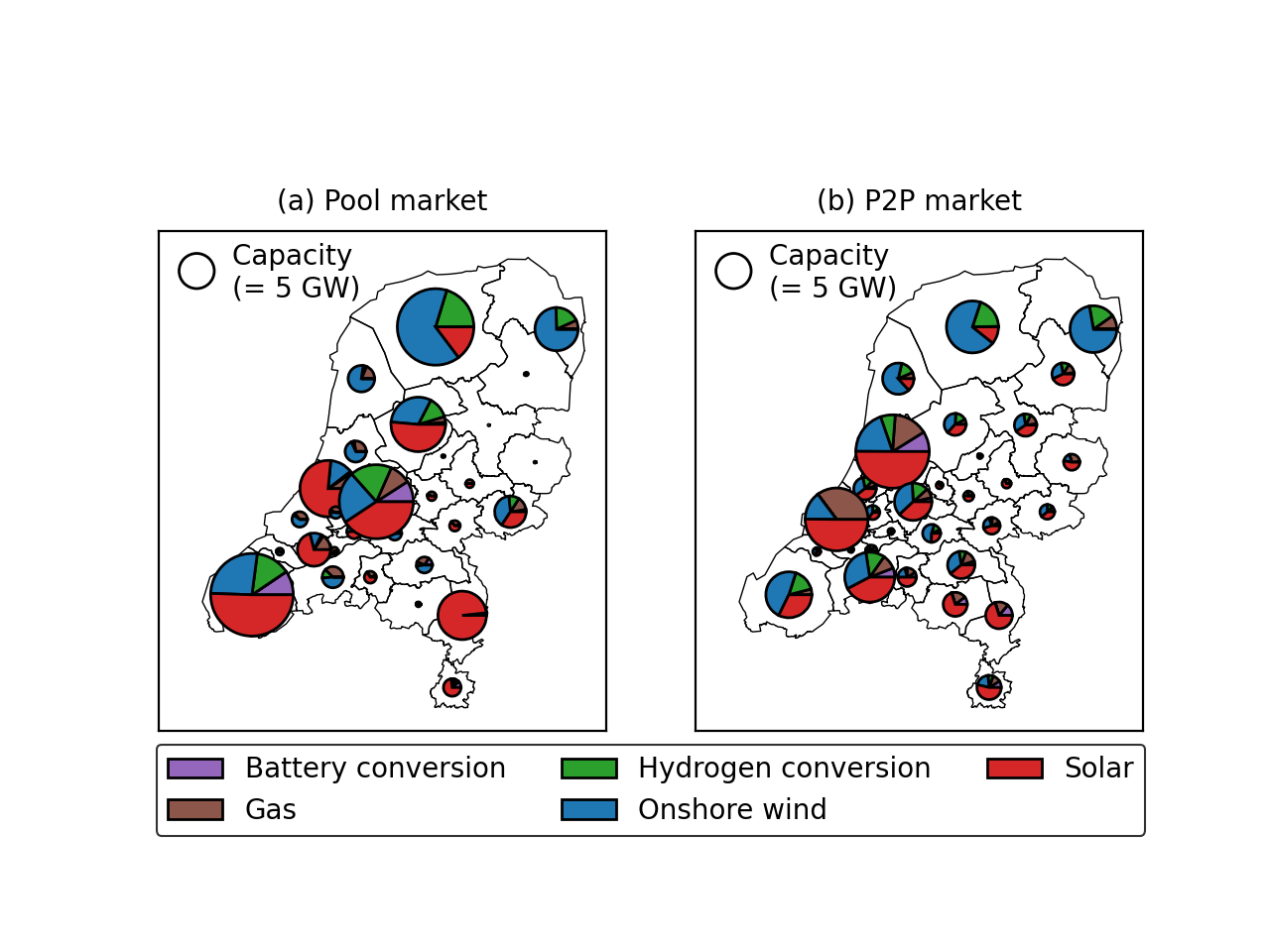}}
\caption{Optimal capacities over the Netherlands. Left to right: pool market, P2P market with 10 \% TC.}
\label{fig:case_I_capacity_map}
\end{figure}

From left to right in Figure \ref{fig:case_I_total_capacity}, the following general trends can be derived. First of all, the capacity of onshore wind declines while the size of solar and storage capacity increases. Secondly, hydrogen storage and battery storage play an equally important role and both have a significant surge in capacity. Thirdly, the expansion of transmission network capacity is marginal for all the cases. To be more specific, the capacity of onshore wind needs to climb to around 35 GW. Even though offshore wind energy is not considered in this model, this indicates how much offshore wind capacity is needed. Given that the capacity factors for offshore wind turbines are generally higher than those of the onshore, the obtained result is more conservative compared to the case when offshore wind is considered. Moreover, to reach the emission goals, the needed capacity for solar PV ranges from 41 GW to 55 GW. Due to the increase in solar capacity, the system levelized cost of energy (LCOE) increases for each case, which are 88 \euro{}/MWh, 90 \euro{}/MWh, 93 \euro{}/MWh, 96 \euro{}/MWh, 97 \euro{}/MWh and 99 \euro{}/MWh, respectively.

Figure \ref{fig:case_I_capacity_map} shows the optimal capacities over the country for the two markets, respectively. 
We start by looking at the planning decisions under the pool market (left figure). Because no preferences are considered in the pool market, cost-optimal results are obtained. Wind capacity is mostly placed on the coastal regions where the wind resources are good, i.e. the northwestern borders. Compared to wind, the capacity factors of solar PV are more evenly spread over the country. 
The right figure shows the planning decisions under the P2P market design with 10 \% TC. The results are significantly different from those in the left figure. Compared to the cost-optimal results, trading barriers are introduced by the use of TC, and thus the resulted capacities are more local, where they are in line with the energy demands of the regions.     

\begin{table*} 
\centering
\caption{Costs for Friesland region.}
\label{tab:market_cost}

\begin{tabular}{c|cc|cc|c}

\toprule
\multirow{2}{*}{Scenarios} & \multicolumn{2}{c|}{Planning phase} & \multicolumn{2}{c|}{Operational phase} & \multirow{2}{*}{Total cost (M\euro{})} \\ 
& Alternatives & Cost (M\euro{}) & Uncertain future markets & Cost (M\euro{}) & \\

\midrule 

Scenario 1 & \multirow{2}{*}{Pool market (alternative 1)} & \multirow{2}{*}{852} & Pool & 108 & 960 \\ 

Scenario 2 & & & P2P & 65 & 917 \\

\midrule 

Scenario 3 & \multirow{2}{*}{P2P market (alternative 2)} & \multirow{2}{*}{589} & Pool & 102 & 691 \\ 

Scenario 4 & & & P2P & 65 & 654 \\

\bottomrule
\end{tabular}
\end{table*}

\begin{figure}[htbp]
\centerline{\includegraphics[width=0.5\textwidth]{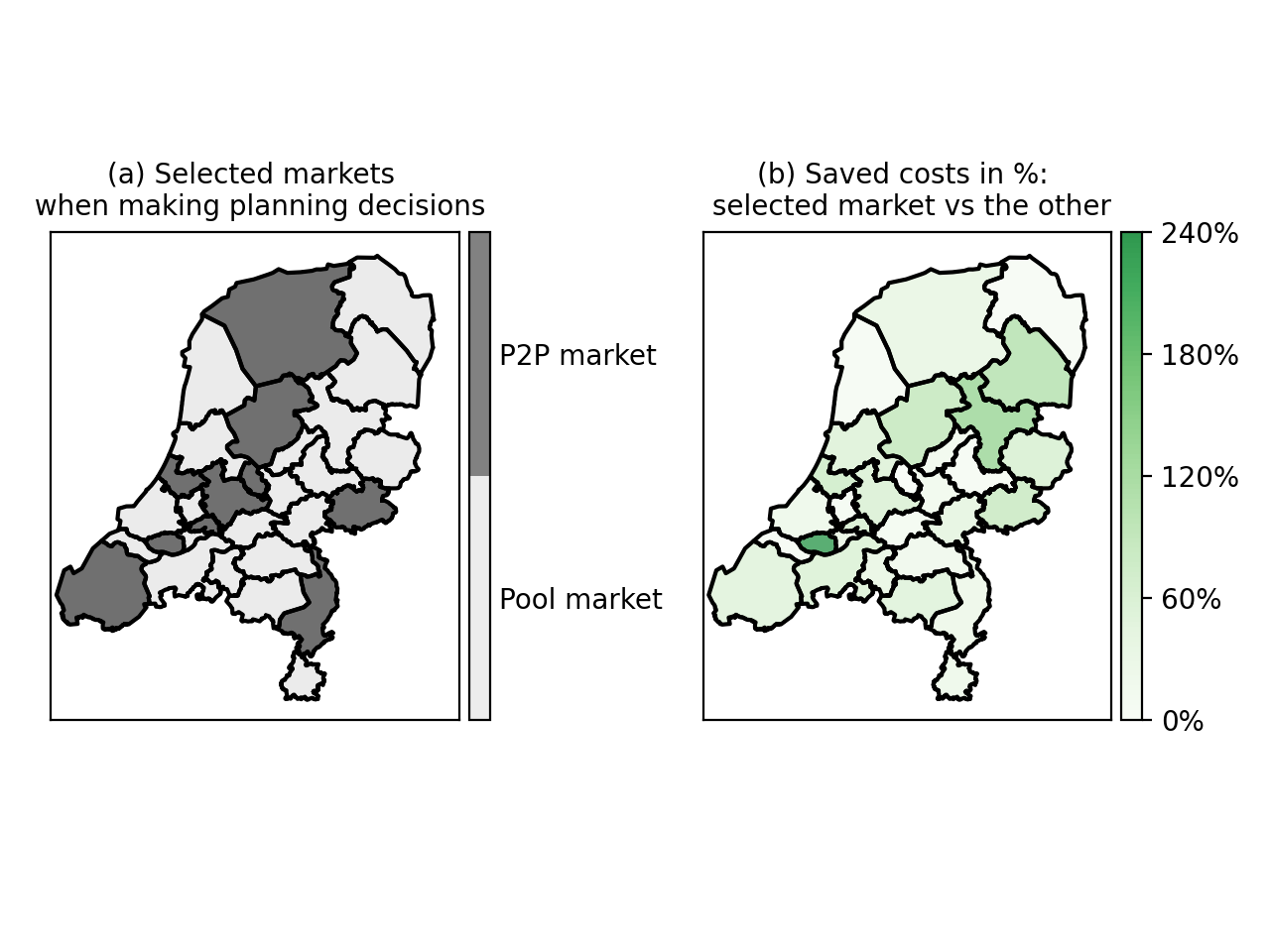}}
\caption{Geographical distribution over the Netherlands of selected market in the planning phase, and the saved costs compared to the other alternative in \%.}
\label{fig:case_I_planning}
\end{figure}

\begin{figure}[htbp]
\centerline{\includegraphics[width=0.5\textwidth]{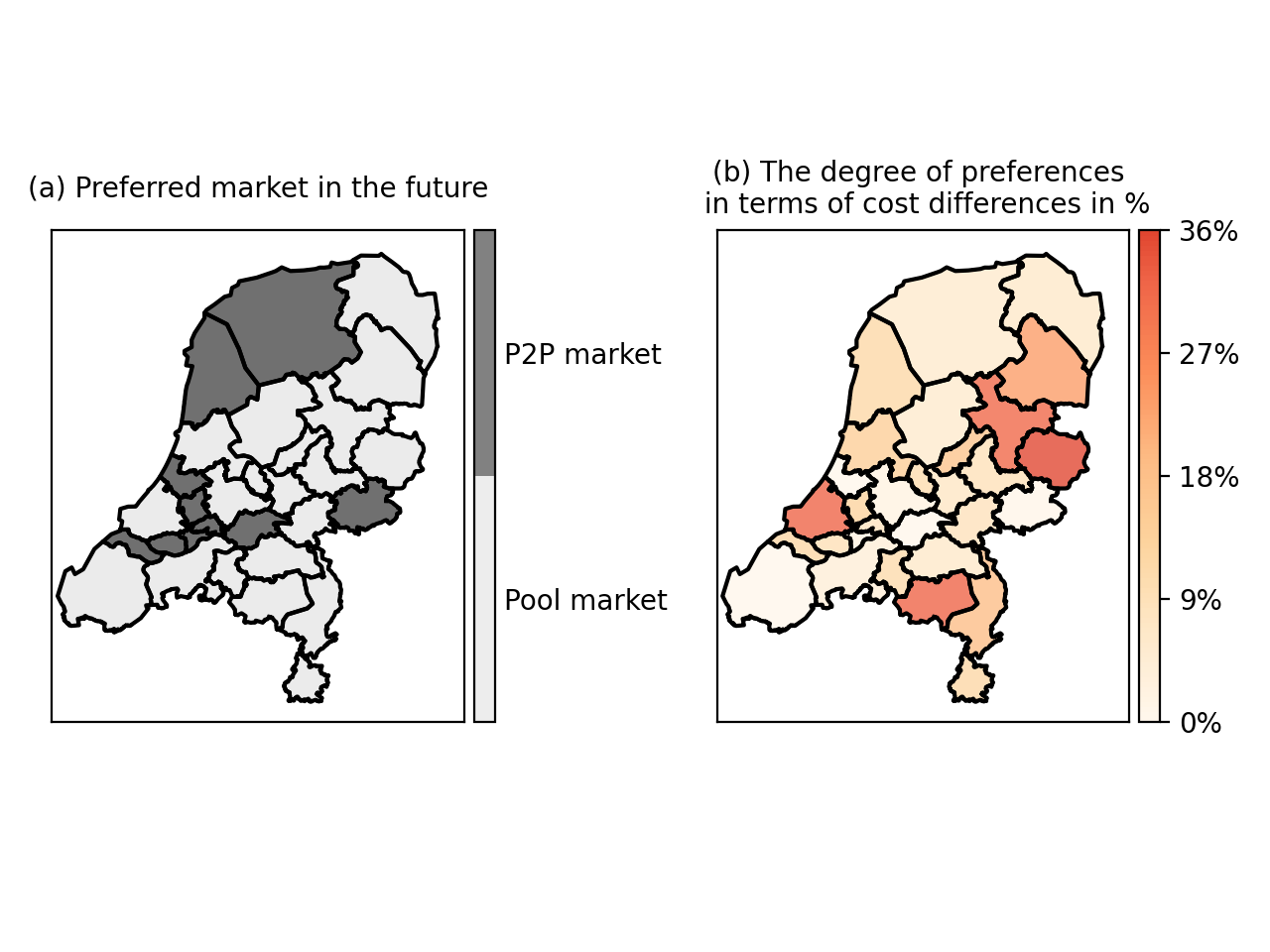}}
\caption{Geographical distribution over the Netherlands of preferred market in the future, and the degree of preferences in terms of the relative cost changes in \%.}
\label{fig:case_I_future}
\end{figure}

\subsection{No-regret planning considering uncertainty in future market designs} \label{sec:case_I_uncertainty}

While it is important to include the trading costs when planning, there is uncertainty in what the future market design will be, especially with the rapid development of the P2P markets. Therefore, it may well happen that planning decisions are made based on a certain type of markets, however, in the future, the market form changes, resulting in cost increases in the operational phase and thus, non-optimal planning decisions. 
In practice, reserve capacities are often used to deal with various uncertainties such as variable demand and intermittent RES supply. In the literature, these uncertainties are often addressed using stochastic or robust optimization methods, whereas the uncertainty in market designs are not studied. In this subsection, we will show how our approach will help to deal with this problem. 

Here, we cope with this issue by using the idea of robust planning. In robust planning, the key idea is that different decisions would lead to different scenarios pertaining to the uncertainties. To assess the scenarios, a min-max problem is solved, meaning that given all the worst-case scenarios (the maximization problem) that may happen in the future, the best decision (the minimization problem) is selected as the solution. In our case, the decision has to be made between two alternatives, i.e. whether to consider the pool market or the P2P market, the uncertainty is the future market design. 

We will now explain the idea further by looking at an exemplary region. The costs of this region is given in Table \ref{tab:market_cost}. 
The discussed region is Friesland, which has good wind resources and low energy demand. 
As explained, one has to do the planning while anticipating a certain market, in this case, either the pool market or the P2P market. The choice of the markets will result in different planning decisions and the anticipated costs. 

As listed in Table \ref{tab:market_cost}, there will be four scenarios. We start with scenario 1 as an example. In this scenario, in the planning phase, the planning decisions are made by assuming a pool market in the future. Then in the future (operational phase), i.e. when the capacities have already been built, it will be the pool market in practice. 
For scenario 1 and scenario 2, the planning decision is made based on the anticipation of a pool market in the future (alternative 1). As a result, many wind capacities have to be built and thus the investment cost is high (852 M\euro{}).
For scenario 3 and scenario 4, the other alternative is taken, which is to anticipate a P2P market in the future, resulting in fewer investments for the region (589 M\euro{}). 

In the operational phase, the electricity market design is uncertain between the pool and the P2P markets. Compared to the P2P market (in scenario 2 and scenario 4), the pool market (in scenario 1 and scenario 3) results in more exports for Friesland, however with lower prices, leading to higher operational costs. 
If these uncertain market designs are not considered, one would only get results of scenario 1 and scenario 4, as we analyzed in the previous sections. However, taking into account the uncertainty, we notice that the costs are different for the four scenarios.
This finding creates difficulties for planning, as it is unclear, which market to consider in order to achieve the lowest total costs.

By listing all the possible scenarios, we are able to find the worst-case scenario when a certain alternative is chosen. If the pool market is anticipated (in scenario 1 and scenario 2), the worst-case scenario is scenario 1 resulting in a total cost of 960 M\euro{}. If the other alternative is chosen, i.e. when the P2P market is considered in the planning phase (in scenario 3 and scenario 4), the worst-case scenario is scenario 3 leading to a total cost of 691 M\euro{}. The solution with the lower cost between the two (691 M\euro{}) is considered as a no-regret solution.
It means that for this region, no matter what the future market design turns out to be, opting for the second alternative (i.e. to plan considering the P2P market) would save 39 \% of the costs compared to the first alternative in the worst case. 

Moreover, even when the second alternative is chosen, there is a cost difference between scenario 3 and scenario 4. This indicates a preference from this region regarding what the future market designs should be. Scenario 4 reveals lower cost compared to that of scenario 3, and thus the P2P market will be preferred, also in the long run, by the region. 
In addition, the absolute value of the cost differences shows the degree of preference. 

At last, we present the results for all the regions. 
Figure \ref{fig:case_I_planning} shows the markets that the regions would consider in the planning phase, and the corresponding saved costs in terms of percentages.
It turns out that a third of the regions would do the planning based on a P2P market, and most of these regions are close to the load centers.
The differences in costs display a wide range from 0 to 240 \%. In particular, the region with the highest cost difference would do the planning based on the P2P market, which, the same as Friesland, is a net export region. 
Figure \ref{fig:case_I_future} illustrates the regions' preferred markets and their degrees of preferences. Although most regions would prefer a pool market in the future, the cost differences are rather small compared to the potential saved costs when choosing alternatives, indicating milder degrees of preferences.

\section{Case Study II} \label{sec:case_II}

\begin{figure*}[htbp]
\centerline{\includegraphics[width=\textwidth]{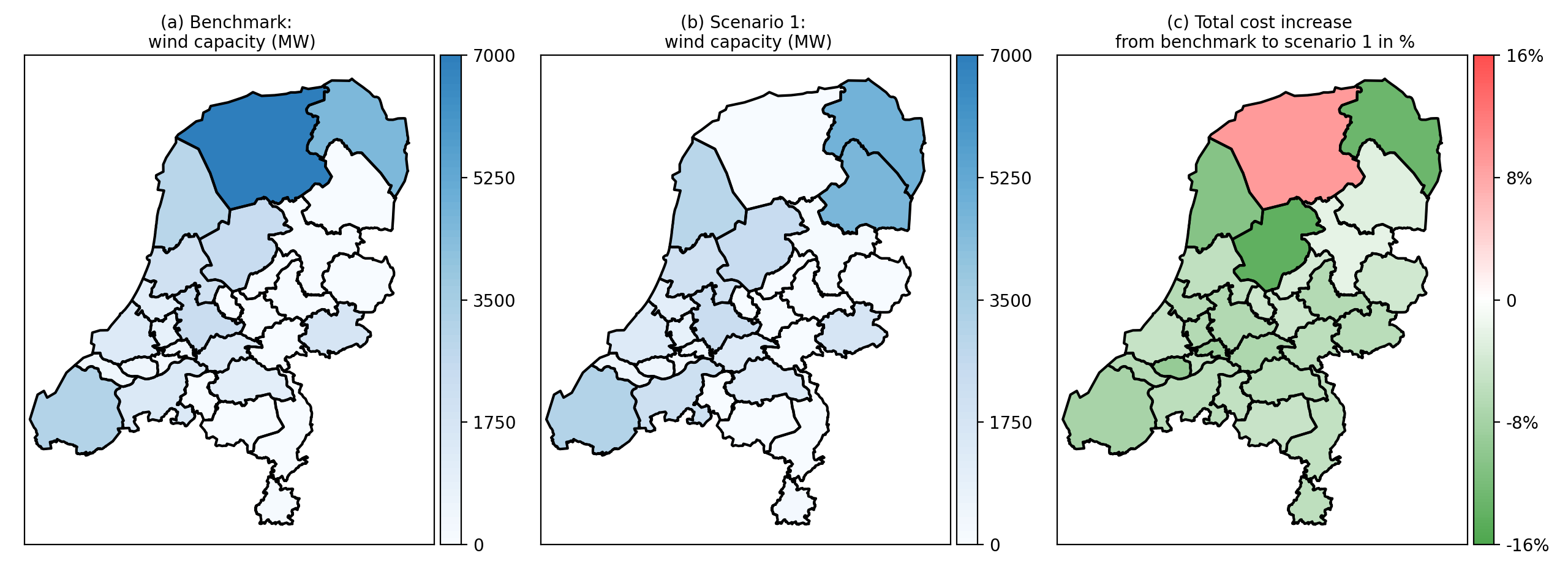}}
\caption{Geographical distributions of wind energy capacity in benchmark and scenario 1, and net increase in total costs in \%.}
\label{fig:case_II_choropleth}
\end{figure*}

In this case study, we again look at the RES planning of the Netherlands, but with a focus on the investment preferences of the regions, in particular related to the wind energy investment.

Due to social acceptance issues, wind energy in the Netherlands is contentious, and thus all the regions are careful with proposing wind energy in their regions. The initial investment proposals show that solar energy is more popular than wind, despite the higher cost. 
Abandoning wind energy will surely make the system LCOE increase which would affect all the regions.
Our framework is capable of taking into account the bilateral relations between the regions, and could be used as a simulator to study the effects of such decisions of not investing in wind energy (results to be shown in Section \ref{sec:case_II_nowind}) and helps with the negotiation process (results to be shown in Section \ref{sec:case_II_counter}).  

\begin{figure}[htbp]
\centerline{\includegraphics[width=0.5\textwidth]{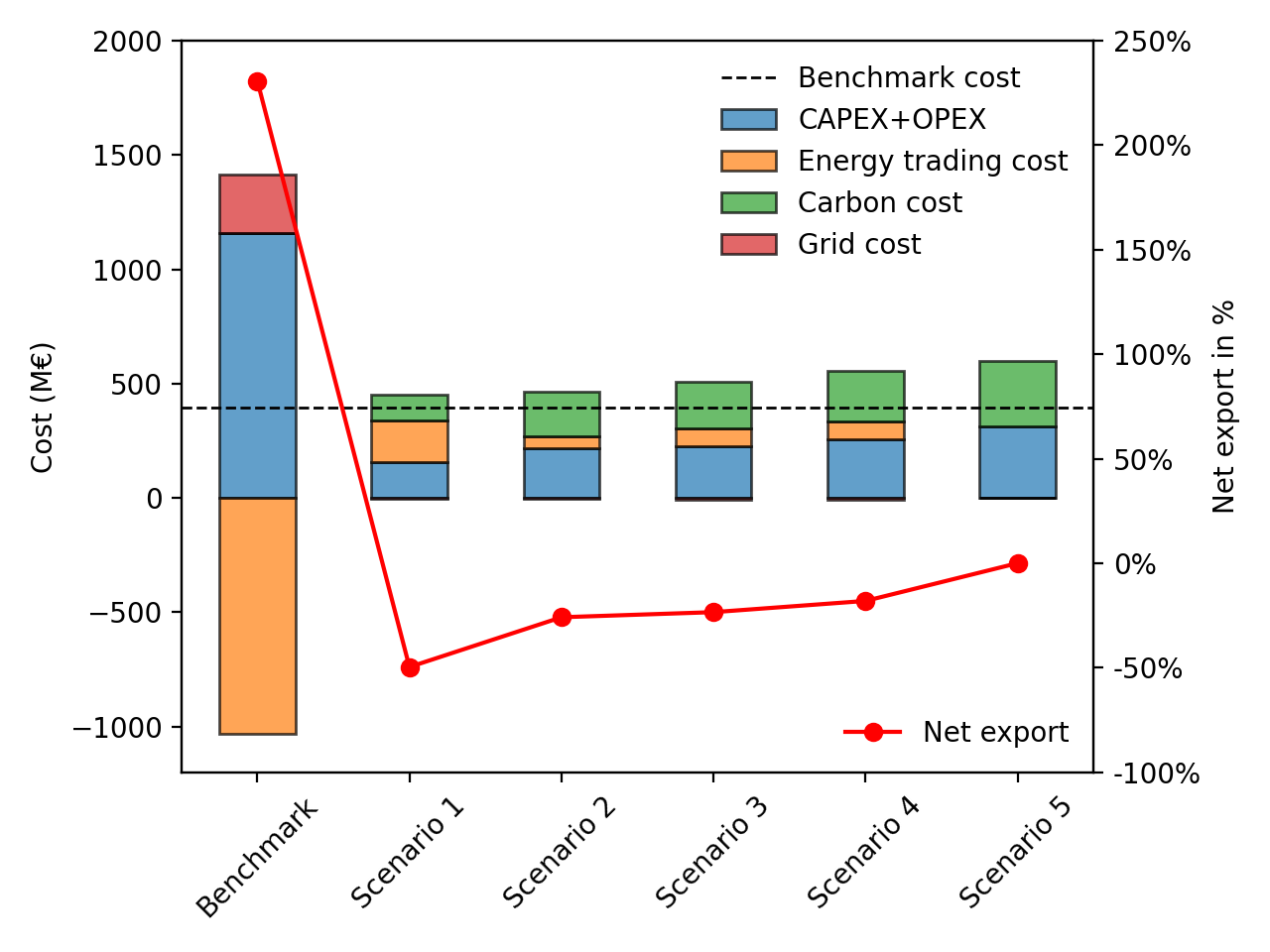}}
\caption{Costs of Friesland in various scenarios.}
\label{fig:case_II_cost_8}
\end{figure}

\subsection{No wind for a region} \label{sec:case_II_nowind}

We first look at the benchmark situation where there is no preferences towards certain RES, i.e. wind energy resistance is not taken into account. Next, the investment preferences against wind energy are considered. These two situations will be referred to as benchmark and scenario 1, respectively.
For illustration purposes, we consider only one region that acts against wind energy. 
The influences of this choice will be shown.

Figure \ref{fig:case_II_choropleth} depicts the wind capacity distribution for the benchmark situation. It shows a concentration in capacity in the coastal regions (later referred to as capacity centers) , i.e. the northwestern border of the country. Due to the discrepancy of load centers and capacity centers, most of the energy that is produced in capacity centers will be transmitted to the load centers which are in the west and middle. 

Among the capacity centers, we look at the Friesland region. It is a farm region and is located in the north. There, energy demands are low (3 \% of the total energy demand) whereas the optimal wind capacity is the largest in the country (20 \% of the total wind capacity). For this reason, as a hypothetical case, we assume that due to high social resistance, this region proposes no investments in wind energy. Other regions or other groups of regions may also be chosen, where the gist of the case study still applies.  
Taking this preference into account, Figure \ref{fig:case_II_choropleth} displays the new capacity distribution. Because of the decrease in wind capacity in Friesland, almost all other regions will have to build more generation capacities. In particular, Drenthe has to build 4680 MW more in wind energy.

Furthermore, we look at the cost changes in percentages. Friesland, as a result of fewer costs in investment yet more costs for importing energy, ends up with a total cost that is 16 \% higher. Although a few of its neighboring regions benefit by gaining more revenues from exporting energy, most of the regions are incurred with more costs due to the large increases in CapEx. This demonstrates that the deliberate choice of one region influences the planning decisions of all other regions. More specifically, most of the regions suffer, though unwillingly, both in terms of increased total costs and forced wind energy investments locally.

\subsection{Other regions' negotiation strategies}\label{sec:case_II_counter}

The RES planning is a negotiation process between the regions. When other regions are unwilling to comply with Friesland, which would be the case in our hypothetical example, their bilateral relationships deteriorate. This, in turn, have an effect on Friesland as well.
Here, our framework is used to simulate the negotiations between the regions. This will be done by considering the preference costs associated with the energy trades which is used as a means to represent the willingness to pay for the prosumers. In this exemplary case, as a result of Friesland's choice, the willingness to pay for other regions with regard to trades with Friesland becomes low. In other words, a high cost is imposed on the trades from other regions as their negotiation strategy to Friesland's proposal.

Figure \ref{fig:case_II_cost_8} shows the costs of Friesland in various scenarios. We start with the benchmark. Due to the large investment needs in wind energy, the CapEx and OpEx costs (together referred to as investment costs later) are high. However, since most of the energy that is generated there will be transmitted and sold to other regions, Friesland will gain a significant revenue from selling energy as well. Overall, its total cost in the benchmark case is 393 M\EUR{}. 
The benchmark case provides the least cost solution for this region. If preferences against wind energy are taken into account, its cost will increase.
In that case, its investment costs decline to 14 \% of those in the benchmark case. Accordingly, due to the lack of local generation capacities, it has to import energy, with the net export percentage dropping from 230 \% to -50 \%. 

\begin{figure*}[htbp]
\centerline{\includegraphics[width=\textwidth]{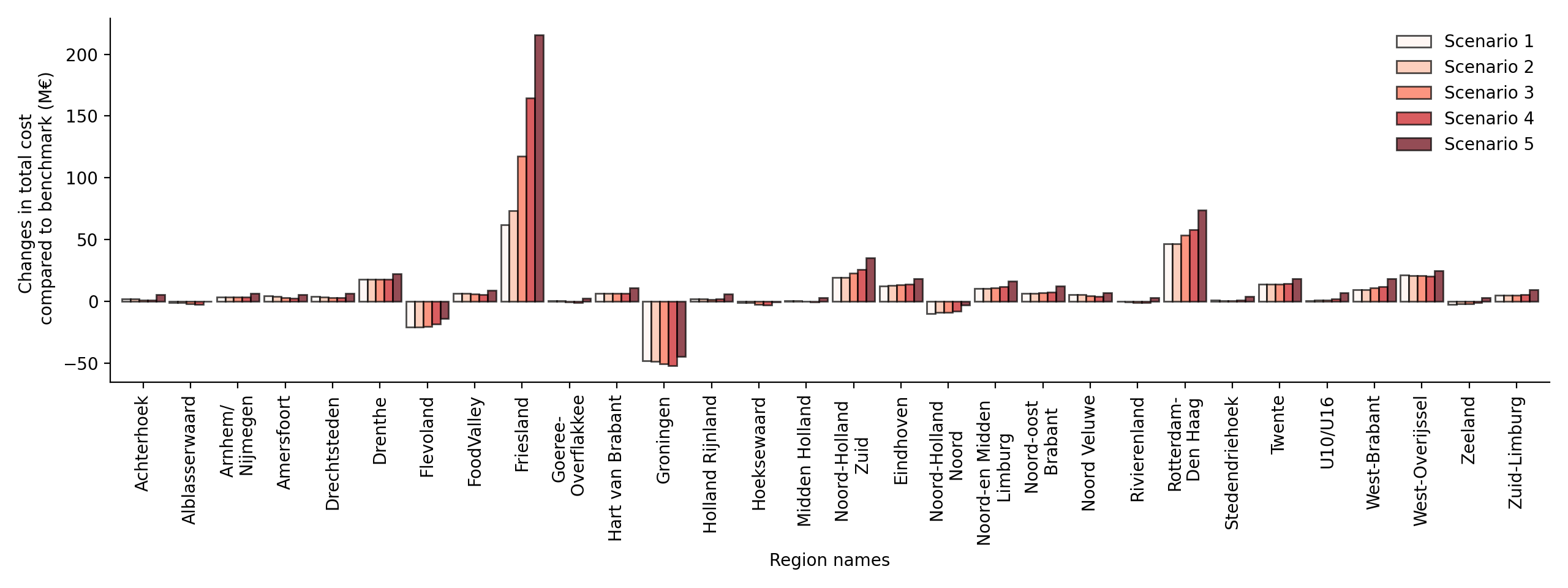}}
\caption{Costs of all the regions compared to benchmark in scenario 1 to scenario 5.}
\label{fig:case_II_cost_all}
\end{figure*}

Then, we look at the results of the following scenarios (2-4) when other regions
start to negotiate with Friesland, with increasing levels of preference cost (10 \EUR{/MWh}, 50 \EUR{/MWh}, 100 \EUR{/MWh}). These indicate degrees of the willingness to pay with trades that involve Friesland. 
The results show that as the trading barriers between other regions and Friesland become larger, the energy trades shrink and thus, Friesland will be forced to rely more on its own energy production which drives up its total costs. In the extreme case (scenario 5), the region will be isolated by others and has no other choice, but to be energy self-sufficient. All these scenarios are not desirable for Friesland, and therefore, it has to reconsider its decision of not investing in wind energy. 

Now we turn to the cost changes for other regions. A key question to answer here is, by imposing trade barriers with Friesland, what are the consequences for other regions? Figure \ref{fig:case_II_cost_all} shows the cost changes relative to the benchmark cost for five cases. There are mainly two groups of regions to be discussed. 
One group is Friesland's neighboring regions with similar wind conditions and low energy demands. Among all, Drenthe builds more wind capacity and is incurred with more costs compared to the benchmark. Flevoland and Groningen have fewer costs, since they benefit from more energy sales. The other group consists of the load centers, Noord-Holland Zuid (Amsterdam region) and Rotterdam-Den Haag, which rely heavily on imports. Due to the choice of Friesland, the energy prices go up leading to higher costs for the load centers as well. In particular, when Friesland does not invest in wind and others take no actions, Rotterdam-Den Haag has a higher cost increase than Friesland. Nevertheless, in all other cases, Friesland bears the most cost increase, especially when the counteraction is strong in negotiations.

\subsection{Discussions}
We would like to emphasize that the topic of this case study is not to propose a national P2P market and study the effects of the P2P markets on the planning decisions nor to argue for the best values to quantify the willingness to pay, but rather to show how the product differentiation term can help to express the preferences of the prosumers and to simulate the inter-regional negotiation process. To this end, we have focused on Friesland as an illustrative region, but we would like to take the discussion away from Friesland into more general inter-regional negotiations. Moreover, we assume all regions have the same willingness to pay, and they all counter one particular region's choices. 
The exact values for willingness to pay depend on the bilateral relationships between other regions and the region under study, which can relate to economic aspects such as how much influence they perceive for their own regions or socio-political aspects such as the political tensions between them. Some regions may even benefit, as shown already. In addition, with various values for willingness to pay from the regions, depending on the results, they may again choose to change their perceptions. Therefore, the actual results highly depend on the case-specific situation when the framework is used in practice. Our case study highlights how our framework can be used to investigate these kinds of policy-relevant challenges.

\section{Case Study III} \label{sec:case_III}

In this case study, we focus on the other part of the framework, i.e. the mixed bilateral and pool market. We expand the geographical context to Europe. Here, the ENTSO-E perspective is taken, where a pan-European electricity market is considered. We follow the model assumptions and techno-economics parameters in \cite{Schlachtberger2017}, where a highly renewable power system was modeled to meet the carbon target in 2050. In that model, 30 countries are connected using HVDC links and a cost-optimal system portfolio was obtained. Compared to the previous case studies, this case study includes more technologies such as offshore wind, hydropower and run-of-river. On top of the cost-optimal solution, to further explore the outcomes of the European power system decarbonization, in this case study, we consider both bilateral and pool markets. In the bilateral market, the countries may choose to trade with others based on their own wills, whereas in the pool market, the energy is traded homogeneously in a pool. Here, bilateral contracts are assumed to account for 70 \% of the energy exchanges of each country \cite{EuropeanCommission2021}. 

\begin{figure*}[htbp]
\centerline{\includegraphics[width=\textwidth]{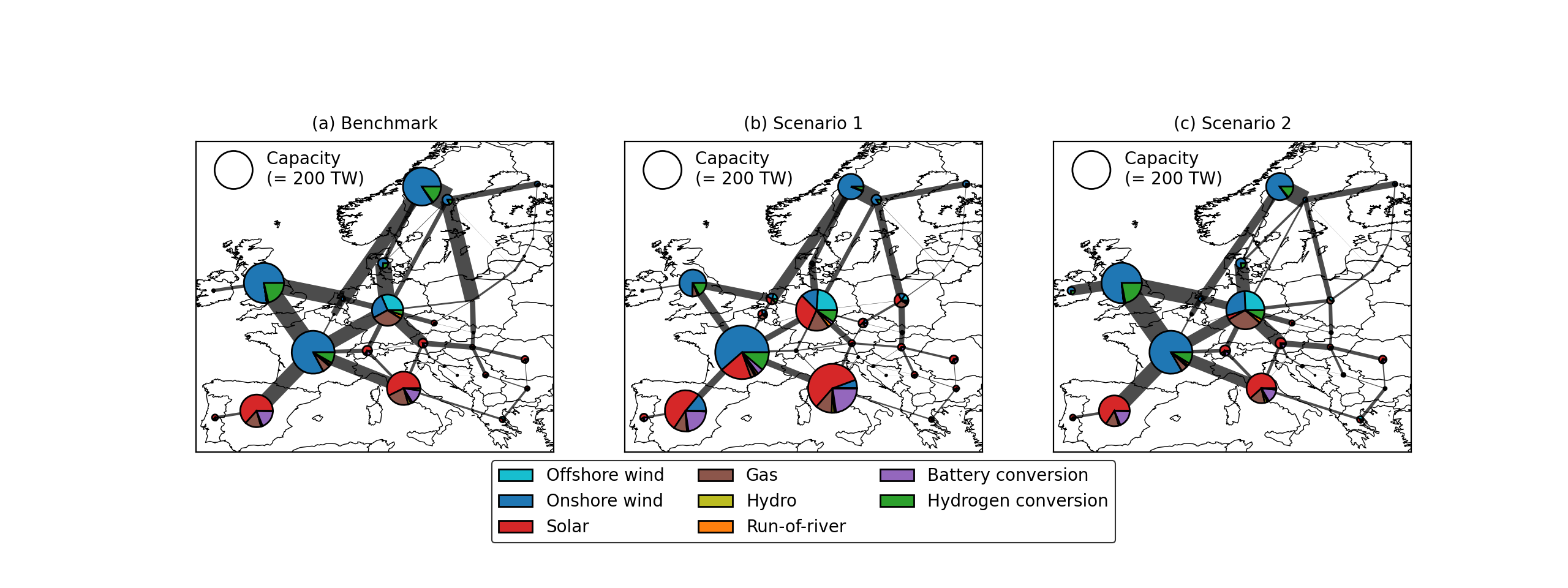}}
\caption{Capacity distribution over Europe for generation, storage and transmission for benchmark, scenario 1 (with non-green index) and scenario 2 (with cross-regional TC). The thickness of the lines represent the transmission capacities.}
\label{fig:case_III_map}
\end{figure*}

\begin{figure}[htbp]
\centerline{\includegraphics[width=0.5\textwidth]{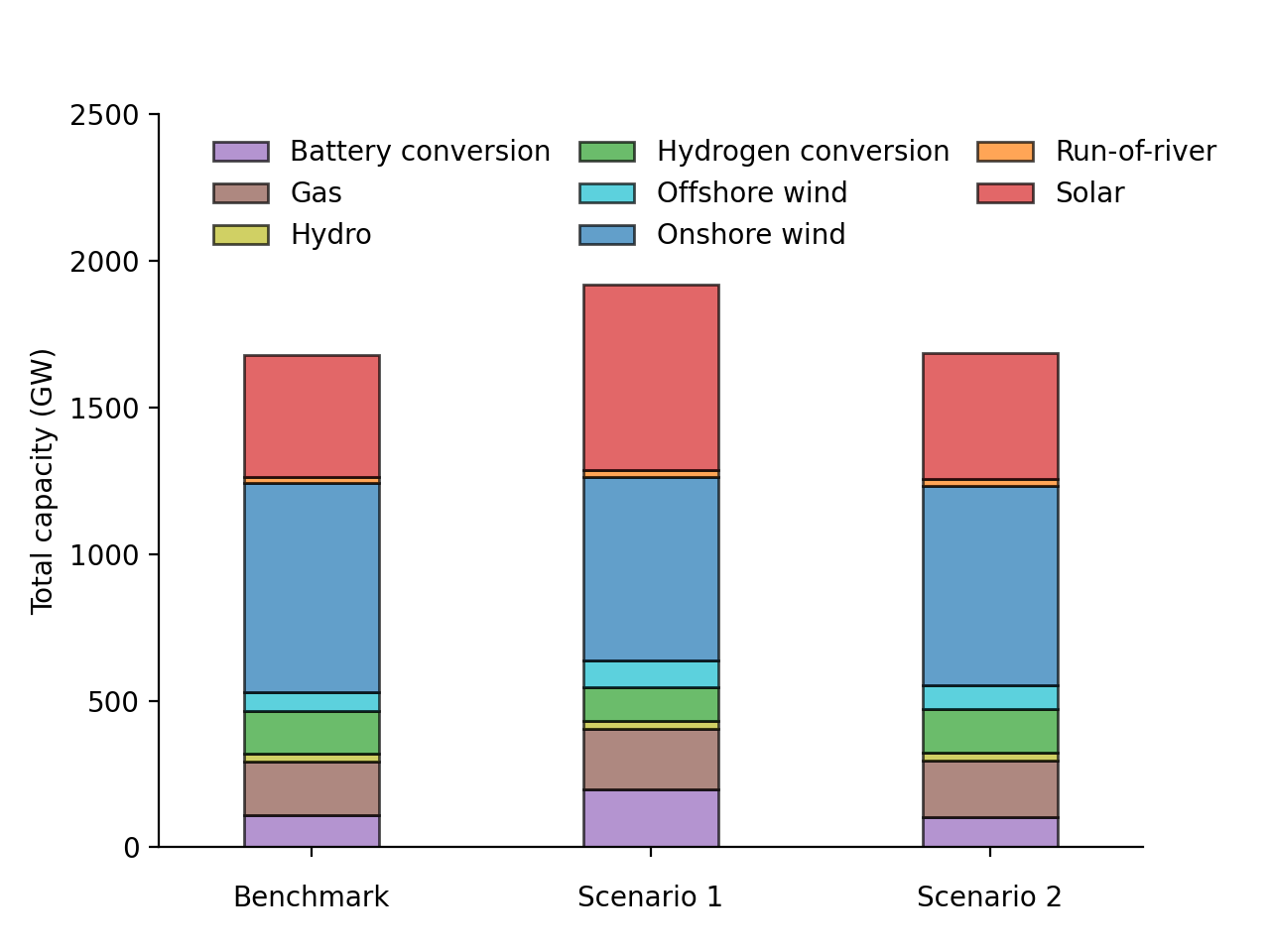}}
\caption{System total capacities for benchmark, scenario 1 (with non-green index) and scenario 2 (with cross-regional TC).}
\label{fig:case_III_total_capacity}
\end{figure}

\begin{figure*}[htbp]
\centerline{\includegraphics[width=\textwidth]{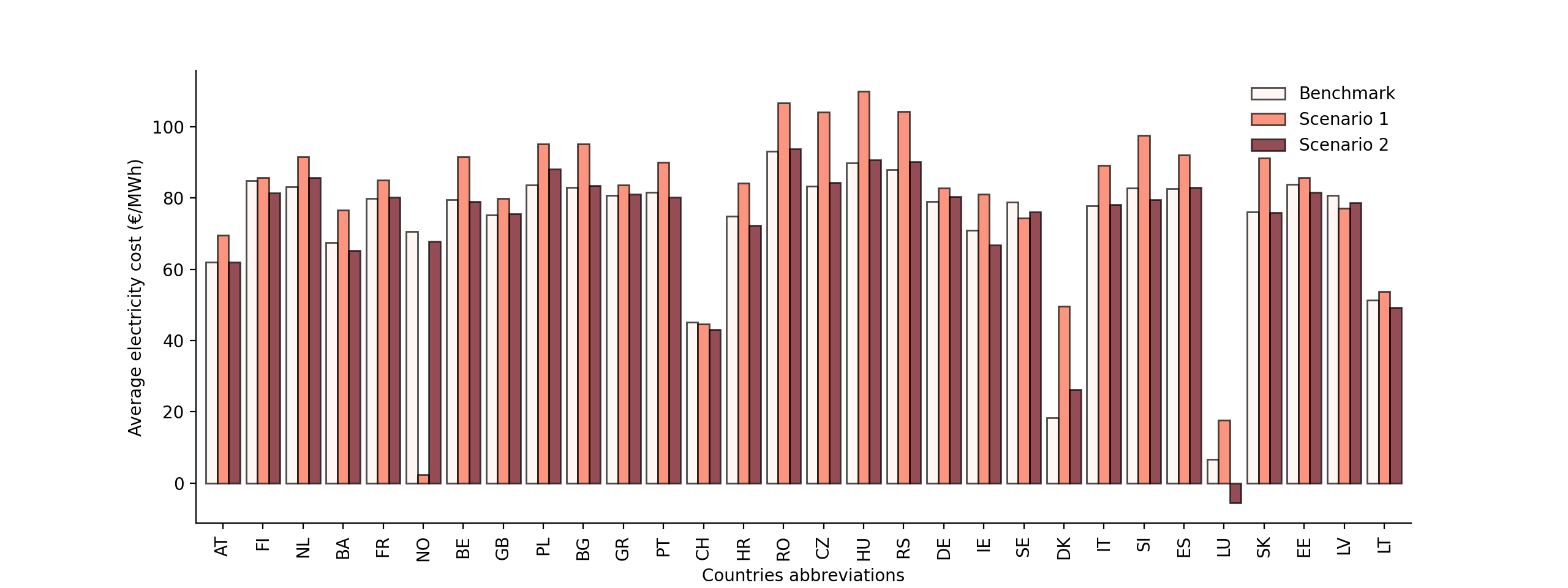}}
\caption{Average costs for the 30 European countries in benchmark, scenario 1 (with non-green index) and scenario 2 (with cross-regional TC).}
\label{fig:case_III_cost}
\end{figure*}

\subsection{Scenario definition}

Three scenarios are investigated in this case study. In addition to a benchmark situation where countries have no preferences in bilateral markets, two other scenarios named scenario 1 and scenario 2 with different usages of the preference terms will be discussed. 

\begin{figure}[htbp]
\centerline{\includegraphics[width=0.5\textwidth]{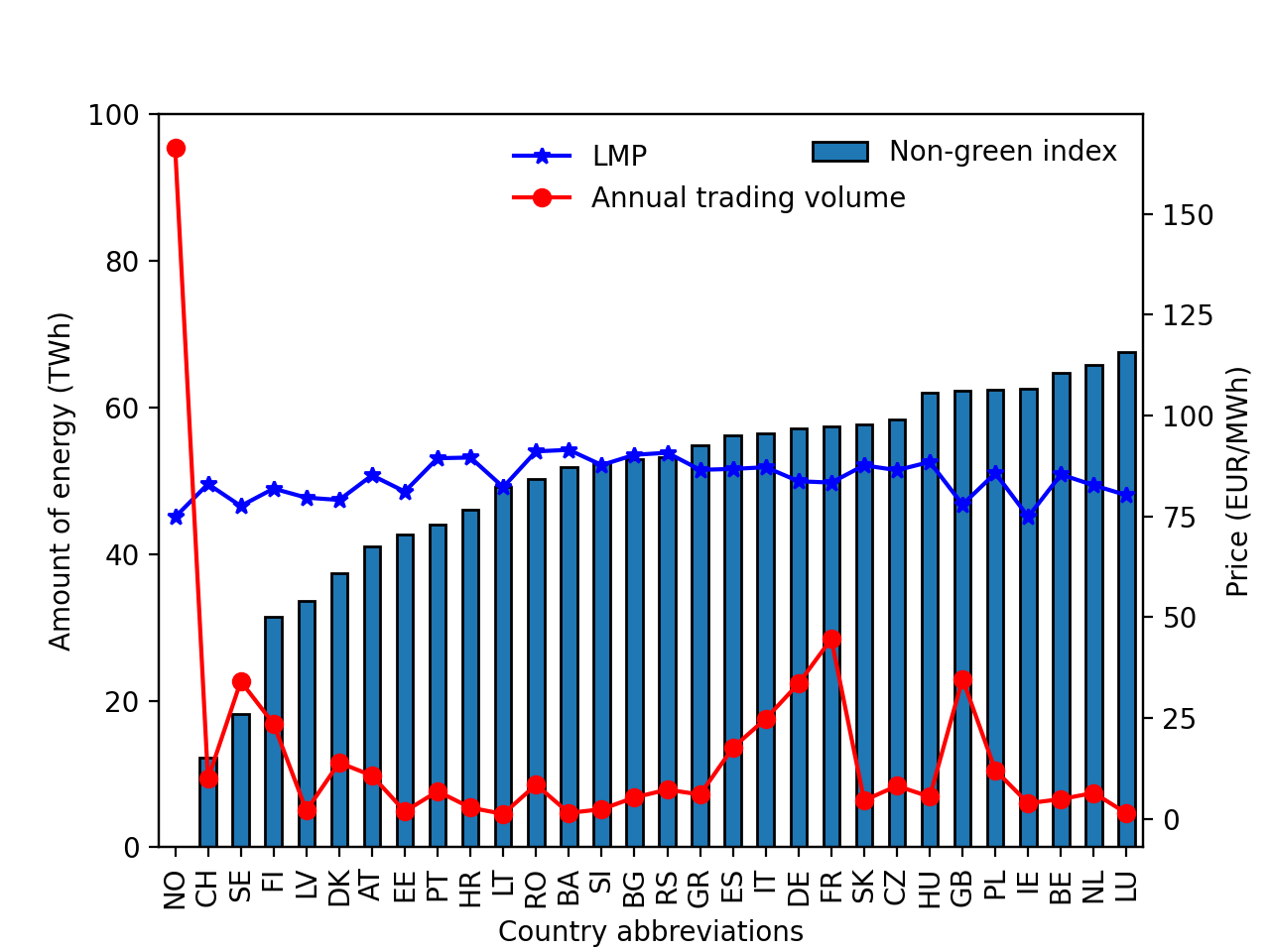}}
\caption{Non-green index, LMP and annual energy trading volumes for each country in scenario 1.}
\label{fig:case_III_price}
\end{figure}

In scenario 1, the preference terms represent bilateral attitudes (or willingness to pay as discussed in Section \ref{case_I_model}) between the countries. This is to show the effects of differentiated preferences on the investment decisions. 
In general, the values of the preference terms are dependent on the geopolitical relationships. We choose to quantify the preferences based on the "non-green index" of each country. 
The idea is that the preference terms are interpreted as a popularity index among the European countries. To be more specific, the countries with a higher percentage of RES in their national generation mix are more popular in bilateral trading for other countries. 
To translate the popularity in costs, we look at the percentage of non-RES generation. For example, Norway is highly dependent on hydropower and only has a non-RES percentage of 25.4 \% in the generation mix. Switzerland ranks second in national RES penetration and has a non-RES share of 37.6 \%. The percentages are translated directly to costs, i.e. 25.4 \EUR{}/MWh and 37.6 \EUR{}/MWh for Norway and Switzerland, respectively. As these costs are derived merely to differentiate the bilateral attitudes, relative values matter more than actual values. We further calibrate all the costs using the lowest one (25.4 \EUR{}/MWh), which ends up with 0 for Norway and 12.2 \EUR{}/MWh for Switzerland. Figure \ref{fig:case_III_price} shows the non-green index of European countries. Luxemburg has the highest preference costs (67.6 \EUR{}/MWh) when trading with others.  

\begin{figure*}[htbp]
\centerline{\includegraphics[width=\textwidth]{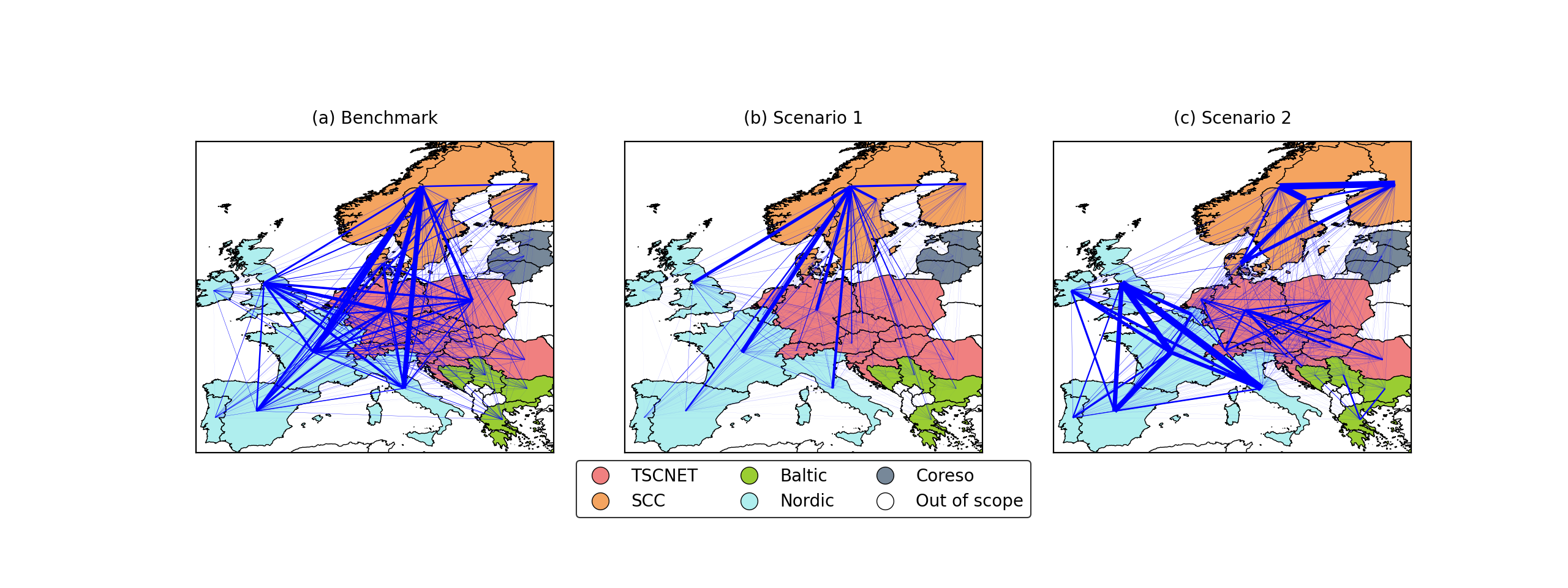}}
\caption{Annual bilateral trading volumes represented by thickness of the lines in benchmark, scenario 1 (with non-green index) and scenario 2 (with cross-regional TC).}
\label{fig:case_III_bilateral}
\end{figure*}

In scenario 2, we use the preference terms to represent the TCs between regions. The TCs could be considered as actual costs or deemed as a measure for cross-regional trading barriers. Different divisions of regions in Europe could be deployed, based on geographical locations or relationships between various multi-national European organizations and agreements. Here, we refer to the TSO regional coordination scheme from ENTSO-E \cite{ENTSO-E2019}. In particular, five regional security coordinators are considered (Figure \ref{fig:case_III_bilateral} shows the geographical locations of the five regions).

\subsection{Investment capacities and costs}
In this subsection, we show the results obtained from the different scenarios. First, the total investment capacity of the system will be given. Then, the geographical distribution of the capacities will be shown. Finally, the average costs of countries are discussed. 

Figure \ref{fig:case_III_total_capacity} illustrates the total investment capacities for the three cases. Several trends could be observed. Firstly, onshore wind energy is the most predominant technology, whose capacity is 42 \%, 33 \% and 40 \% of the generation mix, respectively in the three scenarios. Next, solar energy is the second most needed technology and its capacity evens the onshore wind capacity in scenario 1. 
Then, scenario 3 that encourages intra-regional trades shows a similar generation portfolio with that of the benchmark scenario. Lastly, all other technologies such as offshore wind and storage only contribute a small part to the energy mix.

We now look at the geographical distribution of the capacities over Europe in Figure \ref{fig:case_III_map}. The meteorological conditions of wind and solar energy are drastically different in various parts of Europe. Therefore, in general, wind capacities are built in the middle and the north, whereas solar capacities are concentrated towards the south. Transmission networks are needed to connect the different locations. Most of the generation capacities are concentrated in load centers of Europe, such as France and Germany. Moreover, in the benchmark scenario, it can be found that Norway has built lots of wind energy and thus being a major net export country. In scenario 1 and scenario 2, fewer capacities are built but still many revenues are gained due to the increase in energy price. 

In our framework, the average costs are calculated differently from the traditional LCOE. Normally, when calculating LCOE, only the incurred investment costs are included. However, here, we also include the trading costs, which means that the average costs for a country can even be negative provided that lots of revenues are gained from energy trading. Figure \ref{fig:case_III_cost} depicts the average costs per country for the three scenarios. Overall, the means of the average costs are 73 \EUR{}/MWh, 80 \EUR{}/MWh, 72 \EUR{}/MWh, respectively for the three scenarios. The general trend is also reflected in the costs for most countries, with a few exceptions. For example, Norway is the most popular country to trade with in scenario 1, resulting in the lowest average cost (2 \EUR{}/MWh) therein. Denmark, in contrast, sees a surge in costs in this scenario. This is because it changes from a net exporter to a net importer. It is also found that Luxembourg becomes the only country with a negative average cost in scenario 2. It benefits from the TC-free intra-regional trades and turns its position into a net exporter.

\subsection{Deeper look on the influences of the differentiated non-green index}
After showing the capacities and costs under the three scenarios, we now take a deeper look at how the non-green index influences the results in scenario 1. Figure \ref{fig:case_III_price} illustrates the relationships of the non-green index, the locational marginal price (LMP) and the annual trading volumes for each country. The bar plots show the non-green index representing the preference terms for each bilateral trades. Naturally, one tends to think that the ranking of the non-green index will result in the same ranking for the volumes of the bilateral energy exchanges. This is true for Norway for being the largest exporter, Switzerland's exports are outnumbered by countries such as France, Germany, Sweden and Finland. In fact, the volumes of the bilateral energy exchanges are joint effects of several factors, such as the volumes of local energy production. This is the case with countries where the productions are high, such as France and Germany. Therefore, in terms of export percentages, these countries have lower values than Switzerland. In addition to the national energy production, the other contributing factor is the price of the bilateral trades. The price not only consists of the preference cost, but also the LMP. For the case of Switzerland, its final price might be higher than other countries, which is the other reason to explain its final energy exchange.

\subsection{Further analysis on bilateral trades}
Figure \ref{fig:case_III_bilateral} shows the annual bilateral trades for the three scenarios. In the benchmark scenario, the bilateral exchanges happen actively across Europe, making use of the different wind and solar resources at different locations. In scenario 1, with various preference costs imposed on the trades, the trading volumes drastically decrease with more locally produced energy. Norway, taking advantage of its lowest non-green index, becomes the most active trader across Europe and is selling energy to all other countries. In scenario 2, trading barriers are created by introducing a TC across regions, which notably change the trading trends. Here, bilateral trades are no longer active across regions, but most of them are kept within regions, showing strong intra-regional preferences.

\section{Conclusion} \label{sec:conclusion}

With the rapid development of distributed energy sources, prosumers play an increasingly crucial role in modern power systems. Although peer-to-peer electricity markets have been studied extensively, a largely unrealized value of such markets is to consider their implications on investments. In this paper, we integrate the peer-to-peer markets in planning models, in order to help the prosumers make optimal investment decisions while accounting for their true wishes. Moreover, to systematically assess the influences of market designs on the planning decisions, we consider pool markets and mixed bilateral/pool markets in separate planning models as well.
In Part I of the paper, we develop and elaborate the problem formulations for prosumers, and its coordination with various other actors are also investigated to solve a concurrent generation and transmission planning problem. Here, the roles of the TSO, the energy market operator and the government are revealed. These different yet interconnected optimization problems of the agents are formulated, based on which a centralized problem and a distributed optimization model are derived for each market. 

This framework has been applied to three archetypal case studies. They aim at answering different but all technically, socially and politically relevant questions. For all the cases, the investment decisions and the associated costs are presented. The first case studies the institutional uncertainty in future market designs and helps to obtain no-regret planning decisions for the prosumers. The second case focuses on a situation where a group of prosumers have to decide on their own investments in order to meet a joint carbon target. There, we consider the technology preferences of the prosumers, in particular, an assumed unwillingness to invest in wind energy due to social acceptance issues. The framework acts as a negotiation simulator to inform the group about the consequences of such a preference. Furthermore, by using the bilateral preference terms, it is able to simulate the negotiation strategies that are to the benefit of the prosumers. The third case investigates a European highly renewable power system in 2050. The bilateral relationships between the countries that are shaped by e.g., geopolitical reasons are shown to greatly influence the national planning decisions, and thus the results give practical insights on the decarbonization of the European power system.        

In order to focus the contribution on integrating electricity markets into planning models, some simplifications are made. In order to apply the framework, adjustments could be made from both planning and operational aspects depending on the purpose of the potential applications. These include, e.g., adding more actors such as the distributions system operator, using AC power flow calculation including losses and considering imperfect competition. 
Moreover, the complexity of the distributed optimization algorithms drastically increases, as the number of agents becomes large. It is less of an issue in planning problems, as there are no urgent operational needs to obtain the results within a limited amount of time. Nevertheless, future research could investigate more in this direction to optimize for speed as well as limit information exchange between the agents. In fact, our framework has been constructed in different modules that allow for these modifications.



\bibliographystyle{IEEEtran}
\bibliography{IEEEabrv, ref}

\end{document}